\documentclass[12pt]{article}

\usepackage{amssymb}
\usepackage{amstext}
\usepackage{amscd}

\begin{document}
\newtheorem{teo}{Th\'eor\`eme}[section]
\newtheorem{prop}[teo]{Proposition}
\newtheorem{lem}[teo]{Lemme}
\newtheorem{cor}[teo]{Corollaire}
\newcommand{\OK}{{\cal O}_K}
\newcommand{\IM}{\mbox{Im}}
\newcommand{\Pic}{\mbox{Pic}}
\newcommand{\Spec}{\mbox{Spec}}
\newcommand{\supp}{\mbox{supp}}
\newcommand{\ord}{\mbox{ord}}
\newcommand{\pgcd}{\mbox{pgcd}}
\newcommand{\De}{\Delta_{E/K}}
\newcommand{\oomega}{\omega_{{\cal X}/\OK}}
\newcommand{\Oomega}{\Omega^1_X}

\newcommand{\CCC}{{\mathbb C}}
\newcommand{\RR}{{\mathbb R}}
\newcommand{\ZZ}{{\mathbb Z}}
\newcommand{\QQ}{{\mathbb Q}}
\newcommand{\NN}{{\mathbb N}}
\newcommand{\DD}{{\mathbb D}}
\newcommand{\HH}{{\mathbb H}}

\newcommand{\doo}{\frac{\deg (\omega_{X/O_K})}{[K:\QQ]}}
\newcommand{\nd}{\frac{\log N_{K/\QQ}\Delta_E}{[K:\QQ]}}
\newcommand{\hfi}{\sum_{\sigma} 
\frac{1}{12}\frac{\log \vert\Delta(\tau_{\sigma})\IM(\tau_{\sigma})^6\vert}{[K:\QQ]}}
\newcommand{\OOO}{{\cal O}}
\newcommand{\ocL}{\overline{\cal L}}
\newcommand{\oK}{\overline{K}}
\newcommand{\oL}{\overline{L}}
\newcommand{\oc}{\mbox{\^c}}
\newcommand{\oY}{\overline{Y}}
\newcommand{\oQ}{\overline{\QQ}}
\newcommand{\oCH}{\widehat{CH}}
\newcommand{\oE}{\overline{E}}
\newcommand{\os}{\overline{\sigma}}
\renewcommand{\refname}{R\'ef\'erences}
\renewcommand{\contentsname}{Table des mati\`eres}

\title{Positivit\'e et discretion des points alg\'ebriques des courbes}
\author{Emmanuel Ullmo}
\maketitle

\section{Introduction}

Soient $K$ un corps de nombres et $\oK$ sa cl\^oture alg\'ebrique. Soient $X_K$ 
une courbe propre, lisse, g\'eom\'etriquement connexe de genre $g\ge 2$ sur $K$ et
$J$ sa jacobienne. Soit $D_0$ un diviseur de
degr\'e 1 sur $X$  et $\phi_{D_0}$ le plongement de $X_K$ dans $J$
d\'efini par $D_0$. On note $h_{NT}(x)$ la hauteur de N\'eron-Tate d'un point 
$x\in J(\oK)$. On montre dans ce texte l'\'enonc\'e suivant qui a \'et\'e conjectur\'e
par Bogomolov \cite{Bo}:

\begin{teo}\label{teo1}
Il existe $\epsilon >0$ tel que $\{P\in X_K(\oK) \vert h_{NT}(\phi_{D_O} (P))\le 
\epsilon   \}$ est fini.
\end{teo}

Notons que Raynaud \cite{Ra} a prouv\'e que l'ensemble des
points $P\in X_K(\oK)$ tels que $\phi_{D_0}(P)$ est de torsion dans $J$ est
fini. Le th\'eor\`eme \ref{teo1} g\'en\'eralise cet \'enonc\'e car la condition
$\phi_{D_0}(P)$ est \'equivalente \`a $h_{NT}(\phi_{D_0}(P))=0$. Le lecteur
 s'assurera que la d\'emonstration du th\'eor\`eme \ref{teo1}
est ind\'ependante de celle de Raynaud.

Le th\'eor\`eme \ref{teo1} a \'et\'e obtenu dans de nombreux cas par Szpiro \cite{Sz} 
et Zhang \cite{Za1}. Soit $\Oomega$ le faisceau des diff\'erentielles holomorphes
sur $X_K$. Quand $X_K$ a un mod\`ele propre et lisse sur 
l'anneau des entiers $O_K$ de $K$,
Szpiro \cite{Sz} a montr\'e le th\'eor\`eme \ref{teo1} quand la classe
$[\Oomega-(2g-2)D_0]$ n'est pas de torsion dans $J$ ou la self intersection
$(\omega_{Ar},\omega_{Ar})_{Ar}$ du dualisant relatif au sens d'Arakelov
est non nulle. Il a aussi expliqu\'e comment un \'equivalent du
th\'eor\`eme de Nakai et Moishezon en th\'eorie d'Arakelov permet de prouver
que la non nullit\'e de $(\omega_{Ar},\omega_{Ar})_{Ar}$ est \'equivalente
\`a l'\'enonc\'e du th\'eor\`eme \ref{teo1}. Zhang \cite{Za3} a montr\'e l'analogue
du th\'eor\`eme de Nakai et Moishezon (voir aussi les travaux de Kim \cite{Kim} pour
une autre approche). Notons enfin que toujours dans le cas o\`u $X_K$ admet un
mod\`ele propre et lisse sur $O_K$, Burnol \cite{Bu} et Zhang \cite{Za2} ont 
donn\'e des conditions suffisantes pour que $(\omega_{Ar},\omega_{Ar})_{Ar}$ soit
strictement positif.

  Pour g\'en\'eraliser les travaux de Szpiro
(concernant le cas o\`u le mod\`ele minimal non singulier $\cal X$ de $X_K$ sur
l'anneau des entiers $O_K$ de $K$ est lisse), Zhang \cite{Za1}
 a introduit un accouplement
admissible $(\ ,\ )_a$   g\'en\'eralisant celui d'Arakelov $(\ ,\ )_{Ar}$ sur $\cal X$.
Il a d\'efini un faisceau dualisant relatif $\omega_a$ qui co\"\i ncide avec
le faisceau dualisant relatif $\omega_{Ar}$ dans le cas o\`u $\cal X$ est lisse et 
qui v\'erifie
\begin{equation}
(\omega_{Ar},\omega_{Ar})_{Ar}\ge
(\omega_a,\omega_a)_a \ge 0.
\end{equation}
Il a montr\'e le th\'eor\`eme \ref{teo1} quand $(\omega_a,\omega_a)_a>0$ et quand la
classe $[\Oomega-(2g-2)D_0]$ n'est pas de torsion dans $J$. Il a enfin montr\'e,
 que si $(\omega_a,\omega_a)_a = 0$ et 
$D_0=\frac{\Oomega}{2g-2}$ le th\'eor\`eme \ref{teo1} est en d\'efaut. On obtient
ainsi dans ce texte la preuve du r\'esultat suivant :

\begin{teo}\label{teo4}
Soit $\cal X\longrightarrow \mbox{Spec}(O_K)$ le mod\`ele minimal r\'egulier
d'une courbe lisse g\'eom\'etriquement connexe $X_K$ sur $K$ de genre 
$g\ge 2$. Si ${\cal X}$ a r\'eduction semi--stable alors : 
\begin{equation}
(\omega_{Ar},\omega_{Ar})_{Ar}\ge
(\omega_a,\omega_a)_a>0.
\end{equation}
\end{teo}

Notons que l'in\'egalit\'e $(\omega_{Ar},\omega_{Ar})_{Ar}>0$ a \'et\'e
d\'emontr\'ee par Zhang \cite{Za1} dans le cas o\`u ${\cal X}$ n'est pas
lisse sur $O_K$. Un \'enonc\'e un peu moins g\'en\'eral a \'et\'e obtenu
par Burnol \cite{Bu}. Zhang \cite{Za2} a aussi montr\'e l'in\'egalit\'e  
$(\omega_a,\omega_a)_a>0$ (et donc le th\'eor\`eme \ref{teo1}) quand 
$\mbox{End}(J)\otimes \RR$ n'est pas isomorphe \`a $\RR$, $\CCC$ o\`u
$\HH$ (alg\`ebre des quaternions).

Pour prouver les th\'eor\`emes \ref{teo1} et \ref{teo4}, on suppose que
$(\omega_a,\omega_a)_a = 0$ et que $D_0=\frac{\Oomega}{2g-2}$. On dispose
alors d'une suite infinie de points $x_n$ de $X_K(\oK)$  telle que $h_{NT}
(\phi_{D_0}(x_n))$ tend vers $0$ quand $n$ tend vers l'infini. En utilisant
des th\'eor\`emes d'\'equidistributions des petits points \` a une suite 
$y_n$ de $X_K^g$, construite \`a partir de $x_n$, et \`a son image $z_n$ dans $J$
on obtient une contradiction.

\section{Pr\'eliminaire en th\'eorie d'Arakelov}

\subsection{Vari\'et\'es arithm\'etiques, hauteurs }

\bigskip

Soit $K$ un corps de nombres, $O_K$ son anneau d'entiers et $S=\Spec (O_K)$.
On note $S_{\infty, K}$ l'ensemble des places \`a l'infini de $K$. Une vari\'et\'e
arithm\'etique sur $S$ est la donn\'ee d'un $S$--sch\`ema plat et projectif $X$
dont la fibre g\'en\'erique $X_K$ est lisse.
Pour toute extension
$K_1$ de $K$, on note $X_{K_1}=X_K\otimes_K K_1$. 
 Les points de $X(\CCC)$ vus comme 
$\ZZ$--sch\'ema s'\'ecrivent comme la r\'eunion disjointe $X(\CCC)=\displaystyle
\cup_{\sigma\in S_{\infty,K}} X_{\sigma}(\CCC)$, o\`u $X_{\sigma}(\CCC)=
X\otimes_{\sigma}\CCC$. Un fibr\'e inversible $\oL=(L, \Vert\ \Vert_{\sigma})$
sur $X$ est la donn\'ee d'un fibr\'e inversible $L$ sur $X$ et pour tout
$\sigma\in S_{\infty,K}$ d'une m\'etrique $C^{\infty}$, invariante par la conjugaison 
complexe, sur le fibr\'e inversible $L_{\sigma}=L\otimes_{\sigma}\CCC$ de
$X_{\sigma}(\CCC)$. On note alors $\oL_{\sigma}$ le fibr\'e inversible hermitien
$(L_{\sigma},\Vert\ \Vert_{\sigma})$ de $X_{\sigma}(\CCC)$. On note $\overline{\Pic}(X)$
la cat\'egorie des fibr\'es inversibles hermitiens sur $X$. On dit que 
deux \'el\'ements ${\cal \oL}$, ${\cal \oL}'$ de
$\overline{\Pic}(X)\otimes \QQ$ coincident sur la fibre g\'en\'erique 
s'il existe sur la fibre g\'en\'erique un isomorphisme de ${\cal L}_K$
sur ${\cal L}'_K$ qui est une isom\'etrie
en toute place \`a l'infini. Dans cette situation, on se permet d'\'ecrire
${\cal \oL}_K={\cal \oL}'_K$

En particulier, un fibr\'e inversible hermitien $\oL$ sur $S$ est la donn\'ee
d'un $O_K$--module projectif de rang 1 et de m\'etriques hermitiennes
sur le $\CCC$--espace vectoriel de dimension 1, $L_{\sigma}$, pour toute
place \`a l'infini $\sigma$ de $K$. Le degr\'e d'un fibr\'e inversible hermitien
$\oL$ sur $S$ est d\'efini par l'\'egalit\'e :
\begin{equation}
\deg_{Ar}(\oL)
=\log \#(L/O_K .s)-\sum_{\sigma\in S_{\infty,K}} \log \Vert s\Vert_{\sigma}
\end{equation}
pour une section arbitraire $s$ de $L$.

Soit $X$ une vari\'et\'e arithm\'etique sur $S$ de dimension (absolue) $d$ et
$\oL$ un fibr\'e inversible hermitien sur $X$. Pour tout
$x\in X_K(\oK)$, on note $D_x$ la cl\^oture de Zariski de $x$ dans 
$X$ et $K(x)$ son corps de rationalit\'e  . La hauteur $h_{\oL}(x)$
 est alors d\'efinie par la formule
$\displaystyle h_{\oL}(x)=\frac{\deg_{Ar}(\oL\vert D_x)}{[K(x):K]}$. 
Si on part d'une vari\'et\'e projective lisse $X_K$ sur $K$, que l'on fixe
une extension $K_1$ de $K$, un mod\`ele ${\cal X}_1$ de $X_{K_1}$ sur l'anneau
des entiers $O_{K_1}$ de $K_1$ et un fibr\'e inversible hermitien ${\cal \oL}$
sur ${\cal X}_1$, on peut encore d\'efinir une hauteur sur $X_K(\oK)$. En effet
pour tout point $x\in X_K(\oK)$, $\mbox{Spec}(K(x)\otimes_K K_1)$ est une r\'eunion 
de points $(x_1,\dots,x_r)$ de $X_L(\oL)$. On pose alors 
$$
h_{{\cal \oL}} (x)=\frac{\displaystyle\sum_{i=1}^r h_{{\cal \oL}}(x_i)}{[K_1:K]}.
$$

Soit $T$ une vari\'et\'e analytique complexe de dimension $d-1$ et $\oL=(L,\Vert\ \Vert)$
un fibr\'e inversible hermitien sur $T$. Soit $K$ la forme de courbure
associ\'ee \`a $\oL$ \cite{GH}. On notera dans la suite $c_1(\oL)$ la 
$(1,1)$--forme ferm\'ee $\frac{i}{2\pi}K$. Par ailleurs on notera
$c_1(L)$ la premi\`ere classe de Chern d'un fibr\'e inversible $L$ sur
une vari\'et\'e alg\'ebrique $T$. A travers l'application degr\'e, on identifiera
$c_1(L)^{d-1}$ \`a un entier naturel.

Rappelons  que Gillet et Soul\'e \cite{GS1} \cite{GS3} ont d\'efini pour une
vari\'et\'e arithm\'etique $X$ de dimension $d$ des groupes de Chow
arithm\'etiques $\oCH ^i(X)$ pour tout entier naturel $i$. Quand $X$ est irr\'eductible,
on a $\oCH^0(X)\simeq \ZZ$. On dispose d'une application degr\'e
$$
\oCH^d(X) \longrightarrow \RR.
$$
 Pour tout fibr\'e inversible hermitien $\oL$ sur $X$, on dispose d'une
premi\`ere classe de Chern arithm\'etique $\oc_1(\oL)\in \oCH^1(X)$. Gr\^ace
au produit d'intersection 
$$
\oCH^i(X)\times \oCH^j(X)\longrightarrow
\oCH^{i+j}(X)\otimes_{\ZZ} \QQ,
$$
on sait d\'efinir $\oc_1(\oL)^i\in \oCH^i(X)\otimes_{\ZZ}\QQ$ pour tout
$i\in \NN^*$. On d\'efinit $\oc_0(\oL)=1$ et on voit $\oc_1(\oL)^d$ comme
un nombre r\'eel \`a travers l'application degr\'e.

Rappelons aussi que dans le cas des surfaces arithm\'etiques, on dispose d'une
th\'eorie due \`a Arakelov \cite{Ar}, Faltings  \cite{Fa} et Zhang \cite{Za1} qui 
pr\'ecise les choix de m\'etriques sur les fibr\'es
que l'on est amen\'e \`a \'etudier.
 Soient donc $X$ une courbe lisse et g\'eom\'etriquement connexe 
de genre $g$ non nul sur $K$
et ${\cal X}$ son mod\`ele r\'egulier minimal. Quitte \`a \'elargir $K$, on suppose
que ${\cal X}$ est semi--stable. Pour tout plongement
$\sigma$ de
$K$ dans
$\Bbb C$, on note $X_{\sigma}$ la surface de Riemann obtenue \`a partir de $X$ par le
changement de base d\'efini par $\sigma$. La surface $X_{\sigma}$ est munie d'une
$(1,1)$--forme canonique 
\[
\nu_{\sigma}=\frac{i}{2g}\sum_{i=1}^{g}\omega_i\wedge \overline{\omega}_i,
\]
pour une base orthonorm\'ee $(\omega_1,\dots ,\omega_g)$ de $H^0(X_{\sigma},\Omega^1)$
 o\`u $\Omega^1$ est le faisceau des diff\'erentielles holomorphes sur $X_{\sigma}$ 
pour le produit scalaire~:
\begin{equation}\label{prodscal}
\langle \alpha,\beta \rangle=\frac{i}{2}\int_{X_\sigma}\alpha\wedge \overline{\beta}.
\end{equation}
Arakelov a d\'efini une th\'eorie des intersections $(\ ,\ )_{Ar}$ pour les \'el\'ements
de la cat\'egorie  $\Pic_{Ar}({\cal(X)})$ des
fibr\'es inversibles sur ${\cal X}$ munis en chaque place \`a l'infini $\sigma$ de $K$
d'une m\'etrique permise (\`a courbure proportionnelle \`a $\nu_{\sigma}$).
 Le faisceau
$\oomega$, dualisant relatif de ${\cal X}$ sur $\mbox{Spec} (O_K)$
 est canoniquement muni de
m\'etriques permises \cite{Ar}. On note $\omega_{Ar}=\overline{\oomega}$ l'\'el\'ement 
de $\Pic_{Ar}({\cal(X)})$ ainsi obtenu.

Zhang \cite{Za1} a \'etendu et g\'en\'eralis\'e l'intersection d'Arakelov.
Il a d\'efini une notion d'admissibilit\'e en toute place de $K$ qui co\"\i ncide
avec celle d'Arakelov aux places \`a l'infini. On note 
$\Pic_a({\cal X})$ la cat\'egorie des fibr\'es inversibles admissibles
au sens de Zhang \cite{Za1}. On dispose d'une th\'eorie des intersections                                                              
 $(\ ,\ )_a$ sur $\Pic_a({\cal X})$. 
  Zhang a aussi d\'efini des
m\'etriques admissibles en toute place de $K$ sur le fibr\'e $\oomega$.
On note dans ce texte $\omega_a$ l'\'el\'ement de $\Pic_a({\cal X})$ ainsi
obtenu. On dispose ainsi d'une hauteur $h_{\omega_a}$ sur ${\cal X}$ qui est un 
repr\'esentant de la classe des hauteurs de Weil associ\'es au fibr\'e $\Oomega$
(voir \cite{Si} pour une introdution au hauteurs). On aura besoin du r\'esultat
suivant qui r\'esulte imm\'ediatement de  \cite{Za1} (th\'eor\`eme 2-4). 

\begin{lem}\label{aprox} 
Soit ${\cal X}\longrightarrow \mbox{Spec}(O_K)$ le mod\`ele minimal non singulier
d'une courbe lisse, g\'eom\'etriquement connexe, de genre $g\ge 2$ sur $K$. On suppose
que ${\cal X}$ a r\'eduction semi--stable.
 Il existe une suite d'extension $K_n$ de $K$, 
telle que si on note ${\cal X}_n$ le mod\`ele minimal non singulier
de $X_{K_n}$ sur l'anneau des entiers $O_{K_n}$ de $K_n$,
il existe une suite d'\'el\'ements 
 ${\cal \oL}_n$ de  ${\rm Pic}_{Ar} ({\cal X}_n)\otimes_{\ZZ} \QQ$ telle 
 que pour tout $n\in \NN$,  ${\cal \oL}_n$ co\"\i ncide sur la fibre g\'en\'erique
avec $\omega_{Ar}$ comme fibr\'e inversible hermitien et telle  que 
$$
\sup_{x\in X_K(\oK)} \vert h_{{\cal \oL}_n}(x)-h_{\omega_a}(x)\vert
$$ 
tende vers $O$ quand $n$ tend vers l'infini.
\end{lem}

Dans la suite de ce texte, on fixe une surface arithm\'etique ${\cal X}\rightarrow
\mbox{Spec}(O_K)$ telle que $\omega_a^2=(\omega_a,\omega_a)_a=0$. On choisit un diviseur
$D_0$ sur $X_K$ tel que $D_0=\frac{\Oomega}{2g-2}$. 
Comme $D_0$ 
est fix\'e jusqu'\`a la fin de ce texte on se permet de noter $j=\phi_{D_0}$ 
le plongement de $X_K$ dans sa jacobienne d\'efini par $D_0$ et 
en faisant une extension convenable, on suppose que $D_0$ est 
rationnel sur $K$. On note encore $D_0$ le diviseur horizontal de $\cal X$
de fibre g\'en\'erique $D_0$. Le lemme suivant
r\'esulte imm\'ediatement de \cite{Za1} (preuve du th\'eor\`eme 5-6).

\begin{lem}\label{NT-ARAK}
Pour tout point $P\in X(\oK)$ on  a :
\begin{equation}  
h_{NT}(j(P))=\frac{g}{2g-2}h_{\omega_a}(P).
\end{equation}
\end{lem}

\subsection{Th\'eor\`emes d'\'equidistribution}

Une suite de points $(u_n)$ d'une vari\'et\'e alg\'ebrique, irr\'eductible,
 $X$ est dite 
g\'en\'erique, si $u_n$ converge, au sens de la topologie de Zariski,
 vers le point g\'en\'erique de $X$ (autrement dit,
si pour toute sous-vari\'et\'e stricte $Y$ de $X$, il existe au plus un
nombre fini d'indices $i\in \NN$ tels que $x_i$ soit un point de $Y$).

On a montr\'e dans \cite{SUZ} le th\'eor\`eme d'\'equidistribution des petits points
des vari\'et\'es ab\'eliennes
suivant:

\medskip

\begin{teo}\label{teo2}
Soit $A$ une vari\'et\'e ab\'elienne sur un corps de nombres $K$. Soit $(x_n)$ 
une suite g\'en\'erique 
de points  de $A$ telle que $h_{NT}(x_n)$ converge vers  $0$.
 Soit $O(x_n)$ l'orbite  de $x_n$ sous l'action
du groupe de Galois $G_K=\mbox{Gal}(\oK/K)$. 
Pour toute place \`a l'infini $\sigma$, la suite
$$
\frac{1}{\#O(x_n)}\sum_{x\in \sigma(O(x_n))} \delta_x
$$
converge faiblement vers  
la mesure de Haar de masse totale $1$
 $d\mu_{\sigma}$ de $A_{\sigma}(\CCC)\simeq A\otimes_{\sigma}\CCC $. 
\end{teo}

L'\'enonc\'e suivant qui g\'en\'eralise le th\'eor\`eme
d'\'equidistribution des petits points des vari\'et\'es arithm\'etiques
d\'emontr\'e dans \cite{SUZ} nous sera utile dans la suite de ce texte.

\medskip

\begin{teo}\label{equi}
Soit $X\rightarrow {\rm Spec}(O_K)$
 une vari\'et\'e arithm\'etique de dimension $d$.  Soit $\oL$
 un fibr\'e inversible hermitien 
sur $X$ tel que
$L_K$ soit ample et
$c_1(\oL_{\sigma})$ soit positif pour toute place \`a l'infini
 $\sigma$ de $K$. Soit $h$ une
hauteur de Weil sur $X_K$ associ\'ee \`a $L_K$  
 telle que $h(P)\ge 0$ pour tout $P\in X_K(\oK)$.
On suppose qu'il existe une suite $K_n$ d'extensions finies de $K$, une suite 
${\cal X}_n$ de mod\`eles projectifs de $X_{K_n}$ sur l'anneau des entiers
$O_{K_n}$ de $K_n$ et une suite
 $\oL_n$
d'\'el\'ements de $\overline{{\rm Pic}}({\cal X}_n)\otimes\QQ$,
 telle que pour tout $n\in\NN$ et
pour toute place \`a
l'infini $\sigma$ de $K_n$ on ait:
\begin{equation}
   \oL_n\otimes_{O_{K_n}} K_n=\oL_{K_n}.
\end{equation}
\begin{equation}
\sup_{x\in X_K(\oK)} \vert h_{\oL_n}(x)-h(x)\vert \longrightarrow 0\mbox{ quand }
n\rightarrow \infty.
\end{equation}
Soit $(x_n)$ une suite g\'en\'erique de points
de $X_K(\oK)$ tel que   $h(x_i)$ converge vers $0$.
 Alors pour toute place \`a l'infini $\sigma_0$ de $K$ et toute fonction continue
$f$ sur $X_{\sigma_0}(\CCC)$ la suite
$$
\frac{1}{\# O(x_n)} \sum_{x_n^g\in O(x_n)} f(\sigma_0(x_n^g))
$$ converge
vers 
$\displaystyle \int_{X_{\sigma_0}(\CCC)} f(x) d\mu(x)$ o\`u 
$$
d\mu =\frac{c_1 (\oL_{\sigma_0})^{d-1}}{c_1(L_{\sigma_0})^{d-1} }
 $$
 est vu comme une mesure sur 
$X_{\sigma_0}(\CCC)$ de volume 1.
\end{teo}

{\it Preuve.} La preuve donn\'ee ici est
une simple adaptation de celles des th\'eor\`emes similaires de \cite{SUZ}. 
 Soit $f$ une fonction continue sur $X_{\sigma_0}(\CCC)$ telle que
$$
u_n=\frac{1}{\# O(x_n)} \sum_{x_n^g\in O(x_n)} f(\sigma_0(x_n^g))
$$
 ne converge
 pas vers 
$\displaystyle \int_{X_{\sigma_0}(\CCC)} f(x) d\mu(x)$.

On peut   supposer que :

\medskip

\par a) $\displaystyle \int_{X_{\sigma_0}(\CCC)} f(x) d\mu(x)=0$ (changer $f$
en $f-\displaystyle \int_{X_{\sigma_0}(\CCC)} f(x) d\mu(x)$).

\par b) La suite 
$\displaystyle\frac{1}{\# O(x_n)} \sum_{x_n^g\in O(x_n)} f(\sigma_0(x_n^g))$ 
converge vers une
constante $C<0$ (extraire une sous--suite convergente et
changer si n\'ec\'essaire $f$ en $-f$).

\par c) La fonction $f$ est $C^{\infty}$ sur $X_{\sigma_0}(\CCC)$.
\medskip

  Pour tout r\'eel positif $\lambda$,
 on note $\overline{\OOO}_n(\lambda f)$ le fibr\'e inversible hermitien
$\OOO_{{\cal X}_n}$ de ${\cal X}_n$ 
 muni en toute place \`a l'infini
ne divisant pas $\sigma_0$ de la m\'etrique triviale et en toute place \`a l'infini
 $\sigma$ divisant
$\sigma_0$ de la m\'etrique v\'erifiant $\Vert 1 \Vert_{\sigma}(x)=\exp(-\lambda f(x))$ 
(noter que cela a un sens car pour tout $\sigma$ divisant
$\sigma_0$ on a  $X_{\sigma}\simeq X_{\sigma_0})$. 
On note $\oL_n(\lambda f)=\oL_n \otimes \overline{\OOO}_n(\lambda f)$.
Pour $\lambda$ suffisament petit et pour toute place \`a l'infini $\sigma$ de $K_n$,
 $c_1(\oL_n(\lambda f)_{\sigma})$
est positive (et cela ind\'ependament de $n$.)
On remarque que l'on a pour tout $x\in X_K(\oK)$ :
\begin{equation}
h_{\oL_n(\lambda f)}(x) =h_{\oL_n}(x)+\frac{\lambda}{\#(O(x))}\sum_{x^g\in O(x)} f(x^g).
\end{equation}

Par ailleurs,

\begin{equation}\label{eq6}
\frac{\oc_1(\ocL_n\otimes \OOO_n(\lambda f))^d}{[K_n:K]}=
 \frac{\oc_1(\ocL_n)^d}{[K_n:K]}+\mbox{\rm O}(\lambda^2)
\end{equation}
pour une fonction $\mbox{\rm O}(\lambda^2)$ ind\'ependante de $n$
(utiliser a et le fait que $\ocL_{n,\sigma}$ en tant que fibr\'e inversible
hermitien sur $X_{\sigma}\simeq X_{\sigma_0}$ est ind\'ependant de $n$). 
D'autre part en utilisant
le th\'eor\`eme 5-2 de \cite{Za} et l'existence de la suite $u_n$, on voit que  quand 
$n$ tend vers l'infini, 
$\displaystyle\frac{\oc_1(\ocL_n)^d}{[K_n:K]}$ converge vers $0$.

Soit $\varepsilon$ un nombre r\'eel positif. Par la convergence uniforme
de $h_{\ocL_n}$ vers $h$ et la discussion pr\'ec\'edente, il existe
$N\in \NN$ tel que pour tout $n\ge N$ on a:

\medskip
\begin{equation} 
\displaystyle\sup_{x\in X_K(\oK)}\vert h_{\ocL_n}(x)-h(x)\vert
\leq \varepsilon.
\end{equation} 
\begin{equation}
\displaystyle\vert \frac{\oc_1(\ocL_n)^d}{[K_n:\QQ]d.c_1(L)^{d-1} }
\vert\leq 
\varepsilon.
\end{equation}
\medskip

En utilisant le th\'eor\`eme 5-2 de \cite{Za}, on voit que
\begin{equation}
 \lim_{\i\rightarrow\infty}   (\ h_{\ocL_N}(x_i) +
\frac{\lambda}{\# O(x_i)} \sum_{x_i^g\in O(x_i)} f(\sigma_0(x_i^g))\ )
\ge \frac{(\oc_1(\ocL_N)+\oc_1(\overline{\OOO}(\lambda f))^d}{[K_N:K]d.c_1(L_{K_N})^{d-1}}. 
\end{equation}

On en d\'eduit donc que :

\begin{equation}
\lambda C
\ge O(\lambda^2)-2\varepsilon 
\end{equation}

 En faisant tendre $\varepsilon$, puis $\lambda$ vers $0$, on montre
 que $C\ge 0$. Cette contradiction termine la preuve du th\'eor\`eme \ref{equi}

\section{Suites g\'en\'eriques de petits points de $X^g$}

On rappelle que l'on a fix\'e une surface arithm\'etique semi--stable
${\cal X}\rightarrow \mbox{Spec}(O_K)$, qui est le mod\`ele minimal \r'egulier 
d'une courbe $X_K$ lisse  g\'eom\'etriquement connexe de genre $g\ge 2$ sur $K$,
telle que $(\omega_a,\omega_a)_a=0$.
De plus  $D_0=\frac{\Omega^1_X}{2g-2}$ est  suppos\'e \^etre
 rationnel sur $K$.
En utilisant le corrolaire 5-7 de \cite{Za1} et le lemme \ref{NT-ARAK}
on peut construire une suite g\'en\'erique $(t_k)$ de points de $X_K(\oK)$ telle que
$$
h_{\omega_a}(t_k)=\displaystyle \frac{2g-2}{g}h_{NT}(j(t_k))
$$
 tend vers $0$ quand $k$ tend vers l'infini. On poursuit cette
id\'ee en travaillant sur $X_K^g(\oK)$. On note $s$ l'application de $X_K^g$ dans $J$
telle que : 
$$
s(P_1,\dots,P_g)=j(P_1)+\dots+j(P_g). 
$$

Pour toute extension $L$ de $K$, on note $G_L$ le groupe de Galois de $\oK$ sur 
$L$. Si $K_2$ est une extension galoisienne de $K_1$, on note $\mbox{Gal}(K_2/K_1)$
le groupe de Galois de $K_2$ sur $K_1$.
 Soient $Y$ une vari\'et\'e alg\'ebrique d\'efinie sur $K$ et
$x\in Y(\oK)$. Pour toute extension $L$ de $K$, telle que $L\subset K(x)$, on note
$O_L(x)=G_L.x$ l'orbite sous $G_L$ de $x$. On fait la convention $O_K(x)=O(x)$. Le but
de cette partie est de montrer la proposition suivante.

\begin{prop}\label{suite}
Il existe une suite g\'en\'erique $y_n=(x_{n,1},\dots ,x_{n,g})$ 
de points de $X_K^g(\oK)$
telle que :
\par 1) Pour tout $i\in [1,\dots,g]$, $h_{\omega_a}(x_{n,i})
\longrightarrow 0.$
\par 2) La suite $z_n=s(y_n)$ est une suite g\'en\'erique de $J$ et $h_{NT}(z_n)$
 converge vers 0 quand $n$ tend vers l'infini..
\par 3) L'application $s$ induit une bijection de  $O(y_n)$ sur $O(z_n)$
\end{prop} 

{\it Preuve.} Dans la suite, on fixe un plongement $\sigma_0=id$ de $\oK$ dans $\CCC$ et
on identifie $X_K(\oK)$ \`a un sous--ensemble de 
$X_{\CCC}=X_K\otimes_{\sigma}\CCC$. Pour toute
extension $L$ de $K$, on note $L^c$ la plus petite extension galoisienne de $K$
contenant $L$.     On 
choisit une  distance $d$ sur $X_{\CCC}$ d\'efinissant la topologie complexe.

 Comme cela a \'et\'e indiqu\'e au d\'ebut de cette section, on dispose
d'une suite g\'en\'erique $t_k$
de points de $X_K(\oK)$ telle que $h_{\omega_a}(t_k)$ tend vers 0. Par le th\'eor\`eme
d'\'equidistribution \ref{equi} et
le lemme  \ref{aprox}
 on sait que $\{ O(t_k) \ \vert k\in \NN\} $ est dense
pour la topologie complexe de $X_\CCC $. Pour tout $n\in \NN$, il existe donc un indice
$k\in \NN$ tel que l'on ait \`a la fois :
\begin{equation}\label{x1-1}
h_{\omega_a}(t_k)\le \frac{1}{n}
\end{equation}
\begin{equation}\label{x1-2}
\mbox{ pour tout $x\in X_{\CCC}(\CCC)$ il existe $\alpha\in O(t_k)$ tel que 
$d(x,\alpha) \le \frac{1}{n}$}
\end{equation}

On choisit un tel indice $k$ et on pose $x_{n,1}=t_k$.

\medskip

 On pose alors
$$
\mbox{Gal}({K(x_{n,1})^c}/K)=\{\os_1,\dots,\os_r      \}
$$
et  $X_i=X_K\otimes_{\os_i}K(x_{n,1})^c$. On note encore $\omega_a$
le faisceau dualisant relatif (au sens de Zhang) sur le mod\`ele minimal 
non singulier ${\cal X}_i$ de   $X_{\os_i}$ sur l'anneau des entiers $O_{K(x_{n,1})^c}$
de $K(x_{n,1})^c$.

On dispose par la th\'eorie de Galois
\'el\'ementaire d'un morphisme surjectif:
$$
\pi_k \ :\  \mbox{Gal}({K(x_{n,1},t_k)^c}/K)\ \longrightarrow 
\mbox{Gal}({K(x_{n,1})^c}/K).
$$
Pour tout $i\in [1,\dots,r]$ et tout $k\in\NN$,        
on choisit  $\sigma_{i,k}\in \mbox{Gal}({K(x_{n,1},t_k)^c}/K)$ tel que
$$
\pi_k(\sigma_{i,k})=\os_i.
$$
 On constate que les $\sigma_{i,k}(t_k)$ pour
$i\in [1,\dots,r]$ sont des suites g\'en\'eriques de points de 
$X_i(\oK)$ telles que $h_{\omega_a}(\sigma_{i,k}(t_k))$ tend vers 0.
On en d\'eduit que pour tout $k$ assez grand, pour tout $i\in [1,\dots, r] $ et 
pour tout $x\in X_{\CCC}$ on a :

\begin{equation}\label{x2-1}
h_{\omega_a}(\sigma_{i,k}(t_k)) \le \frac{1}{n}.
\end{equation}

\begin{equation}\label{x2-2}
\mbox{  Il existe
 $ \alpha_i\in O_{K(x_{n,1})^c}(\sigma_{i,k}(t_k))$ tel que }
d(x,\alpha_i)\le \frac{1}{n} 
\end{equation}

\begin{equation}\label{x2-3}
\mbox{dim }H^0(X_{\CCC},{\cal O}(x_{n,1}+t_k))=1
\end{equation}

(Remarquer pour ce dernier point que la suite $t_k$ est g\'en\'erique et que 
$\{ P\in X(\CCC) \ \vert  \  \mbox{dim}\ H^0(X_{\CCC},{\cal O}(x_{n,1}+P)) >1   \}$
est fini). On fixe un tel $k$ et on pose $x_{n,2}=t_k$ et $\sigma_{i,k}=\sigma_i$.
 La suite $(x_{n,1},x_{n,2})$
a la propri\'et\'e suivante :

\begin{lem}
Pour tout $(x,y)\in X_{\CCC}\times X_{\CCC}$, il existe $(\alpha_1,\alpha_2)
\in O(x_{n,1},x_{n,2})$ tel que $\max(d(x,\alpha_1),d(y,\alpha_2))\le \frac{1}{n}.$
\end{lem}
{\it Preuve}. D'apr\`es (\ref{x1-2}), il existe $\os_i\in \mbox{Gal}({K(x_{n,1})^c}/K)$
 tel
que $d(x,\os_i(x_{n,1}))\le \frac{1}{n}$. Par ailleurs d'apr\`es (\ref{x2-2}) il existe 
$\gamma \in  \mbox{Gal}({K(x_{n,1},x_{n,2})^c}/K(x_{n,1})^c)$ tel que 
$$
d(y,\gamma \sigma_i(x_{n,2}))\le \frac{1}{n}.
$$
 On constate que $(\alpha_1,\alpha_2)=
\gamma \sigma_i ((x_{n,1},x_{n,2}))$ convient.

\medskip

En proc\'edant de m\^eme, on construit la suite $y_n=(x_{n,1},\dots,x_{n,g})$ de 
$X_K^g$ telle  que 
$$
\mbox{dim } H^0(X_{\CCC},{\cal O}(x_{n,1}+\dots+x_{n,g}))=1
$$
(utiliser le fait que $t_k$ est g\'en\'erique et \cite{Mi} lemme 5-2),
 $h_{\omega_a}(x_{n,i})$ converge
vers 0 pour tout $i$ et telle que pour tout $(x_1,\dots,x_g)\in X_{\CCC}^g$, il existe
$$
(\alpha_1,\dots,\alpha_g)\in O(x_{n,1},\dots,x_{n,g})
$$
 v\'erifiant
$\displaystyle\max_i( d(x_i,\alpha_i))\le \frac{1}{n}$. Cette derniere propri\'et\'e
nous prouve que la suite $y_n$ est g\'en\'erique. 

La deuxi\`eme partie de la proposition
se d\'eduit de la premi\`ere et du fait que la suite $y_n$ 
est g\'en\'erique en utilisant
le lemme \ref{NT-ARAK} et le fait que $s$ est surjective. 

 Comme $D_0$ est rationnel sur $K$, $s$ induit 
une surjection de $O(y_n)$ sur $O(z_n)$. Soit
 $(\alpha_1,\dots,\alpha_g)\in O(y_n)$ tel que 
$$
z_n=x_{n,1}+\dots+x_{n,g}=\alpha_1+\dots+\alpha_g
$$
On sait que les fibres du morphisme $s$ correspondent aux syst\`emes lin\'eaires
sur  $X$ (voir \cite{Mi} chap\^\i tre 5 dans le cas o\`u $D_0$ est un point
rationnel sur $K$ et se ramener \`a ce cas apr\`es translation). Comme 
$$
\mbox{dim } H^0(X_{\CCC},{\cal O}(x_{n,1}+\dots+x_{n,g}))=1,
$$
 on voit que $s$
est fini au dessus de $z_n$. Pour $n$ assez grand la construction de la suite 
prouve que $x_{n,i}$ ne peut pas \^etre dans $O(x_{n,j})$ pour $i\neq j$. On a 
donc $\alpha_i=x_{n,i}$ pour tout $i$ et $s$ est injective. Ceci termine la
preuve de la derni\`ere partie de la proposition.

\section{Preuve du th\'eor\`eme \ref{teo1}}

Dans la suite, on notera $\pi_i$ la $i$--\`eme projection de ${\cal X}^g=
{\cal X}\times_{O_K}\dots\times_{O_K}{\cal X}$ sur ${\cal X}$. On notera aussi
$\pi_i$ la $i$--\`eme projection de $X_K^g$ sur $X_K$. On rappelle que $\nu$ d\'enote
la $1$-$1$--forme canonique sur $X_{\CCC}$.

\begin{lem}\label{sxg}
Soit $y_n$ la suite pr\'ec\'edemment construite. La suite
$$
\frac{1}{\#(O(y_n))}\sum_{(\alpha_1,\dots,\alpha_g)\in O(y_n)} 
\delta_{(\alpha_1,\dots,\alpha_g)}
$$
converge faiblement vers la mesure  $\pi_1^*\nu\wedge\dots\wedge \pi_g^*\nu$ de
$X_{\CCC}^g$.
\end{lem}

{\it Preuve}. 
 On note $\oL=\pi_1^*\oomega\otimes\dots\otimes\pi_g^*\oomega$  le
 fibr\'e inversible m\'etris\'e sur ${\cal X}^g$
obtenu en  munissant $L_K=\pi_1^*\Oomega\otimes\dots\otimes\pi_g^*\Oomega$, 
en toute place \`a l'infini, de la m\'etrique produit des
images inverses des m\'etriques canoniques sur $\Oomega$. 
On note $h$ le repr\'esentant de la classe des hauteurs Weil, sur $X_K^g(\oK)$
 associ\'e \`a 
$L_K$ v\'erifiant :

$$
h(P_1,\dots,P_g)= \sum_{i=1}^g h_{\omega_a} (P_i)
$$
Pour tout $(P_1,\dots,P_g)\in X_K^g(\oK)$. On a $h(P_1,\dots,P_g)\ge 0$.

On fixe dans la suite une suite d'extension $K_n$ de $K$ et une suite
${\cal \oL}_n$ sur le mod\`ele minimal r\'egulier ${\cal X}_n$
de $X_{K_n}$ donn\'ees par le lemme \ref{aprox}. On pose alors
$$
{\cal X}^n_g ={\cal X}_n\otimes_{O_{K_n}} \dots\otimes_{O_{K_n}} {\cal X}_n
$$
$$
\oL_n=\pi_1^*{\cal \oL}_n\otimes\dots\otimes\pi_g^*{\cal \oL}_n
$$
On constate alors que $\oL_n$ co\"\i ncide g\'en\'eriquement 
(en tant que fibr\'e inversible hermitien) avec $\oL_{K_n}$ et que l'on a pour
tout $(x_1,\dots,x_g)\in X_K(\oK)^g$ :
$$
h_{\oL_n}(x_1,\dots,x_g)=\sum_{i=1}^g h_{{\cal \oL}n}(x_i).
$$
On en d\'eduit alors que

$$
\sup_{x\in X_K^g(\oK)}
\vert h_{\oL_n}(x)-h(x)\vert
$$
tend vers $0$ quand $n$ tend vers l'infini.
Or $L_K$ est ample sur $X_K$,
 en toute place \`a l'infini $\sigma$ de $K$, $c_1(\oL)_{\sigma}$
est une $1$-$1$--forme positive et   
on a construit une suite g\'en\'erique $y_n$ telle que $h(y_n)$ tend
vers 0.  Le th\'eor\`eme d'\'equidistribution \ref{equi} nous
permet alors de conclure que la suite 

$$
\frac{1}{\#(O(y_n))}\sum_{(\alpha_1,\dots,\alpha_g)\in O(y_n)} 
\delta_{(\alpha_1,\dots,\alpha_g)}
$$

converge faiblement vers la mesure 
$$
\frac{c_1(\oL)^g}{c_1(L_K)^g}=
\pi_1^*\nu\wedge\dots\wedge \pi_g^*\nu.
$$

\begin{lem}\label{Haar}
Soit $z_n$ la suite de points de $J$ pr\'ec\'edemment construite. La suite
$$
\frac{1}{\#(O(z_n))}\sum_{z\in O(z_n)} \delta_z
$$
converge faiblement vers la mesure de Haar  normalis\'ee $d\mu_{H}$ de $J$.
\end{lem}

{\it Preuve}. C'est une cons\'equence imm\'ediate de la deuxi\`eme
partie de la proposition \ref{suite}, du th\'eor\`eme 
d'\'equidistribution \ref{teo2} des
petits points des vari\'et\'es ab\'eliennes et du lemme \ref{NT-ARAK}.

\begin{lem}\label{eg-mes}
Pour toute fonction $f$,  continue  sur $J$, \`a valeurs dans $\RR$, on a:
\begin{equation}
\int_J f(z) d\mu_H =\int_{X_{\CCC}^g} f. s(x_1,\dots,x_g) 
\pi_1^*\nu\wedge\dots\wedge \pi_g^*\nu.
\end{equation}
\end{lem}

{\it Preuve}. En utilisant le lemme \ref{Haar}, on voit que la suite
$$
\frac{1}{\#(O(z_n))}\sum_{z\in O(z_n)} f(z)
$$
converge vers $\displaystyle\int_J f(z) d\mu_H$. Or d'apr\`es la troisi\`eme
partie de la proposition \ref{suite} on a :
\begin{equation}
\frac{1}{\#(O(z_n))}\sum_{z\in O(z_n)} f(z)
=\frac{1}{\#(O(y_n))}\sum_{(\alpha_1,\dots,\alpha_g)\in O(y_n)}
 f. s(\alpha_1,\dots,\alpha_g).
\end{equation}
D'apr\`es le lemme \ref{sxg},
cette derni\`ere somme converge vers
$$
\int_{X_{\CCC}^g} f. s(x_1,\dots,x_g) 
\pi_1^*\nu\wedge\dots\wedge \pi_g^*\nu.
$$

\begin{lem}\label{contradiction}
On \`a l'\'egalit\'e :
\begin{equation}
s^*d\mu_H=g!\ \pi_1^*\nu\wedge\dots\wedge \pi_g^*\nu.
\end{equation}
\end{lem}
{\it Preuve}. Cela r\'esulte du lemme \ref{eg-mes} et de la formule
de changement de variables quand on a remarqu\'e que $s$ est une application
g\'en\'eriquement finie de degr\'e $g!$.

{\it Preuve du th\'eor\`eme \ref{teo1}}. 
On va prouver le th\'eor\`eme \ref{teo1} en montrant que
cette derni\`ere \'egalit\'e ne peut \^etre r\'ealis\'ee. Pour cela explicitons
ces deux mesures. Soit $(\omega_1,\dots,\omega_g)$ une base orthonorm\'ee de
$H^0(X_{\CCC},\Oomega)$ pour le produit scalaire d\'efini dans l'\'equation \ref{prodscal}.
On note $\omega^J_i$ l'unique forme diff\'erentielle de type $1$-$0$ sur $J$,
invariante par translation, telle que $j^*\omega^J_i=\omega_i$. On a vu que l'on a :
$$
\nu=\frac{i}{2g}\sum_{i=1}^g \omega_i\wedge\overline{\omega_i}
$$
On pose alors 
$$
\mu=\frac{i}{2g}\sum_{i=1}^g \omega^J_i\wedge\overline{\omega_i}^J.
$$
On a alors $d\mu_H=\frac{g^g}{g!}\mu^g$ et $\nu=j^*\mu$ et pour tout
$i\in \{1,\dots,g \}$ on a :
$$
s^*\omega_i^J=\sum_{j=1}^{g} \pi_j^*\omega_i
$$
On trouve alors
\begin{equation}
s^*d\mu_H=\frac{g^g}{g!}(\sum_{i=1}^g \pi_i^*\omega_1)\wedge
 (\sum_{i=1}^g \pi_i^*\overline{\omega_1})
\wedge\dots\wedge (\sum_{i=1}^g \pi_i^*\omega_g)\wedge 
(\sum_{i=1}^g \pi_i^*\overline{\omega_g}) 
\end{equation}

On voit donc que pour
 tout point $P=(P_1,\dots,P_g)\in X^g$ telle que 
$$
H^0(X_{\CCC},\Oomega(-P_1-\dots-P_g))>0
$$
(par exemple $P=(P_0,\dots,P_0)$ avec $P_0$ point de Weierstrass de $X_{\CCC}$)
on a $s^*d\mu_H(P)=0$ . Par contre la forme de type $g$-$g$,
$$
\pi_1^*\nu\wedge\dots\wedge \pi_g^*\nu
$$
est partout strictement positive. Ceci prouve la contradiction du 
lemme \ref{contradiction} et termine la preuve du th\'eor\`eme.


\begin{thebibliography}{99}


\bibitem{Ar} S. J. Arakelov,
{\em Intersection theory of divisors on an arithmetic surface}, Math. USSR--Izv. 
{\bf 8} (1974), 1167--1180.


\bibitem{Bo} Bogomolov, F.A. {\it Points of finite order on an abelian variety},
Math. USSR Izv. {\bf 17} (1981).


\bibitem{Bu} Burnol, J.-F. {\it Weierstrass points on arithmetic surfaces},
Invent. Math. {\bf 107} (1992) 421--432.



\bibitem{Fa} G. Faltings,
{\em Calculus on arithmetic surfaces}, Ann. of Math.{\bf 119}, 387--424.



\bibitem{GS1} Gillet, H et Soul\'e C.  {\it Characteristic classes for
algebraic vector bundles with hermitian metrics, I, II}, Ann. of Math.
{\bf 131} (1990), 162--203.  

\bibitem{GS3} Gillet, H et Soul\'e C. {\it Arithmetic intersection theory},
Publ. Math. IHES. {\bf 72} (1990) 94--174.


\bibitem{GH} Griffiths, P et Harris, J. {\it Principles of algebraic geometry},
Wiley Interscience, (1978).

\bibitem{Kim} Kim, M. {\it Small points on constant arithmetic surfaces},
Duke. Math. Journal. {\bf 61}, vol {\bf 3} (1990)  823--833. 


\bibitem{Mi} Milne, J.S.  {\it Jacobian Varieties} dans "Arithmetic Geometry",
\'edit\'e par Cornell et Silverman, Springer Verlag (1985), 167--212. 


\bibitem{Ra} Raynaud, M. {\it Courbes sur une vari\'et\'e ab\'elienne 
et points de torsion}, Invent.math {\bf 71} (1983), 207--223.


\bibitem{Si} Silverman, J.H.  {\it The Theory of Height Functions} 
dans "Arithmetic Geometry",
\'edit\'e par Cornell et Silverman, Springer Verlag (1985), 151--166. 



\bibitem{Sz} Szpiro, L. {\it Sur les propri\'et\'es num\'eriques du dualisant
relatif d'une surface arithm\'etique}, The Grothendieck Festchrift. {\bf III}
Progress in Mathematics (1990).

\bibitem{SUZ} Szpiro, L. Ullmo, E. Zhang, S.  {\it Equidistribution des petits points}
\`a para\^\i tre dans Invent. Math.


\bibitem{Za3} Zhang, S. {\it Positive line bundles on arithmetic surfaces},
Ann. of Maths. {\bf 136} (1992), 569--587.


\bibitem{Za1} Zhang, S. {\it Admissible pairing on a curve},
Invent. Math. {\bf 112} (1993) 171--193.

\bibitem{Za} Zhang, S. {\it Positive line bundles on arithmetic varieties},
J. Amer. Math. Soc. {\bf 8} (1995)  187--193.


\bibitem{Za2} Zhang, S.  {\it Small points and adelic metrics},
J. Algebraic Geometry. {\bf 4} (1995) 281--300.



\end{thebibliography}
\end{document}